\documentclass{amsart} 
\usepackage{amssymb, amsmath}
\usepackage[dvips]{graphicx}
\usepackage{amsfonts}
\usepackage{latexsym}
\usepackage{color}
\newtheorem{theorem}{Theorem}	
\newtheorem{lemma}{Lemma}[section]		
\newtheorem{corollary}{Corollary}		
\newtheorem{proposition}{Proposition}		
			
\newtheorem{definition}{Definition}	

\newtheorem{claim}{Cliam}[section]
\title[Strictly convex Wulff shapes and $C^1$ convex integrands]
{
Strictly convex Wulff shapes \\ 
and $C^1$ convex integrands 
}
\author{
Huhe Han}
\address{Graduate School of Environment and Information Sciences,Yokohama National University, {Yokohama 240-8501,} Japan}
\email{han-huhe-bx@ynu.jp}
\author{
Takashi Nishimura
}
\address{
Research Institute of Environment and Information Sciences,  
Yokohama National University, 
Yokohama 240-8501, Japan}
\email{nishimura-takashi-yx@ynu.jp}
\begin{document}
\begin{abstract}
{
In this paper, it is shown that a Wulff shape is strictly convex 
if and only if its convex integrand is of class $C^1$.      
Moreover, applications of this result are given.   

} 
\end{abstract}
\subjclass[2010]{\color{black}52A20, 52A55, 82D25} 
\keywords{Wulff shape, convex integrand, convex body, dual Wulff shape, spherical Wulff shape, 
spherical convex body, spherical dual Wulff shape.  
} 
\maketitle  

\section{Introduction}
\par
\begin{figure}[hbtp]
\begin{center}
\includegraphics[width=6cm]{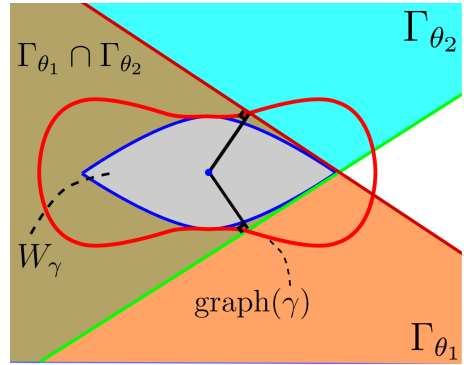}
\caption{A Wulff shape $\mathcal{W}_\gamma$.}
\label{figure 1}
\end{center}
\end{figure}

Let $\gamma: S^n\to \mathbb{R}_+$ be a continuous function, where 
$n$ is a positive integer, $S^n$ is the unit sphere in $\mathbb{R}^{n+1}$ and 
$\mathbb{R}_+=\{a\in \mathbb{R}\; |\; a>0\}$.    
For any $\theta\in S^n$, we set 
\[
\Gamma_{\gamma, \theta}=\{x\in \mathbb{R}^{n+1}\; |\; x\cdot \theta\le \gamma(\theta)\}, 
\]    
where the dot in the center stands for the standard dot product of two vectors 
$x, \theta\in \mathbb{R}^{n+1}$.    
Then, the following set $\mathcal{W}_\gamma$ is called the \textit{Wulff shape} associated with 
$\gamma$ 
(see Figure \ref{figure 1}).  
\[
\mathcal{W}_\gamma=\bigcap_{\theta\in S^n}\Gamma_{\gamma, \theta}.
\]  
\par 
The Wulff shape $\mathcal{W}_\gamma$, which was firstly introduced by G.~Wulff in \cite{wulff}, 
is known as a geometric model of a crystal at equilibrium 
(for instance, see \cite{giga, crystalbook, taylor, taylor2}).  
By definition, the Wulff shape $\mathcal{W}_\gamma$ is compact, convex 
and it contains the origin of 
$\mathbb{R}^{n+1}$ as an interior point; namely,  
$\mathcal{W}_\gamma$ is a convex body such that the origin is contained in int$(\mathcal{W}_\gamma)$, 
where int$(\mathcal{W}_\gamma)$ stands for the set consisting of interior points of $\mathcal{W}_\gamma$    
(for details on convex bodies, see \cite{schneider}).    
Conversely, it is known that any convex body 
containing the origin as an interior point is a Wulff shape associated 
with an appropriate support function (\cite{taylor}).     
Thus, for any convex body $W\subset \mathbb{R}^{n+1}$ such that 
$\mbox{int}(W)$ contains the origin,  
there exists the non-empty set, denoted by $C^0_{W}(S^n, \mathbb{R}_+)$, 
consisting  of continuous functions  
$\gamma: S^n\to \mathbb{R}_+$ such that $\mathcal{W}_\gamma=W$.      
Moreover, it is known that for any convex body $W$ such that 
the origin is contained in $\mbox{int}(W)$,  
there exists the smallest element $\gamma_{{}_W}\in C^0_{\small W}(S^n, \mathbb{R}_+)$ 
in the sense that  
$\gamma_{{}_W}(\theta)\le \gamma(\theta)$ is satisfied for any 
$\theta\in S^n$ and any $\gamma\in C^0_{\small W}(S^n, \mathbb{R}_+)$ (\cite{taylor}).   
The function $\gamma_{{}_W}$ is called  
the \textit{convex integrand} of $W$   
( for details on convex integrand, see Section \ref{section 2}).  
\begin{theorem}\label{theorem 1}
Let $W\subset \mathbb{R}^{n+1}$ be a convex body 
containing the origin of $\mathbb{R}^{n+1}$ as an interior point of $W$.   
Then, $W$ is strictly convex 
if and only if its convex integrand $\gamma_{{}_{W}}$ is of class $C^1$.            
\end{theorem}
A more restricted dual relationship than the one given in 
Theorem \ref{theorem 1} has been obtained by F.~Morgan as follows.   
\begin{theorem}[\cite{morgan}]\label{theorem 2}
Let $W\subset \mathbb{R}^{n+1}$ be a convex body 
containing the origin of $\mathbb{R}^{n+1}$ as an interior point of $W$.   
Then, $W$ is uniformly convex  
if and only if its convex integrand $\gamma_{{}_{W}}$ is of class $C^{1,1}$.            
\end{theorem}
\noindent 
For the definitions of uniform convexity and of class $C^{1,1}$, see 
\cite{morgan}.  
As explained in p. 348  of \cite{morgan}, 
the notion of strict convexity (resp., class $C^1$) is 
certainly weaker than the one of 
uniform convexity (resp., class $C^{1,1}$) for 
Wulff shapes (resp., convex integrands).   
Moreover, \lq\lq strict convexity\rq\rq (resp., \lq\lq class $C^1$\rq\rq ) 
is more common and easy to treat than \lq\lq uniform convexity\rq\rq  
(resp., \lq\lq class $C^{1,1}$\rq\rq).      
Thus, Theorem \ref{theorem 1} may be regarded 
as a useful genelarization of 
Theorem \ref{theorem 2}.       
%
%
\par 
\bigskip 
In Section \ref{section 2}, preliminaries are given.  
Proof of Theorem \ref{theorem 1} is given in Section \ref{section 3}.     
In Section \ref{section 4}, 
applications of Theorem \ref{theorem 1} are given.     
\section{Preliminaries}\label{section 2}
\subsection{Convex integrands}\label{convex integrands}
\par
\begin{figure}[hbtp]
\begin{center}
\includegraphics[width=6cm]{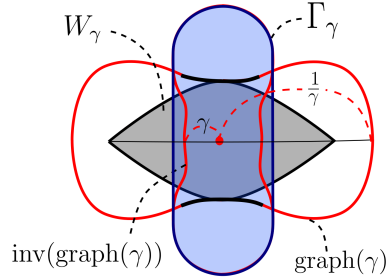}
\caption{A Wulff shape $\mathcal{W}_\gamma$ and the convex hull $\Gamma_\gamma$.}
\label{figure 2}
\end{center}
\end{figure}

 
Let $\gamma: S^n\to \mathbb{R}_+$ be a continuous function.   
Set 
\[
\mbox{graph}(\gamma)=\{(\theta, \gamma(\theta))\in \mathbb{R}^{n+1}-\{0\}\; |\; \theta\in S^n\}, 
\]
where $(\theta, \gamma(\theta))$ is the polar plot expression for a point of $\mathbb{R}^{n+1}-\{0\}$.   
Let $\mbox{inv} : \mathbb{R}^{n+1}-\{0\}\to \mathbb{R}^{n+1}-\{0\}$ be the inversion with respect to 
the origin of $\mathbb{R}^{n+1}$, namely, 
$\mbox{\rm inv}: \mathbb{R}^{n+1}-\{0\}\to \mathbb{R}^{n+1}-\{0\}$ is defined as follows:   
\[
\mbox{inv}(\theta, r)=\left(-\theta, \frac{1}{r}\right). 
\]
Let $\Gamma_\gamma$ be the boundary of the convex hull of 
$\mbox{inv}(\mbox{graph}(\gamma))$ 
(see Figure \ref{figure 2}).
\begin{definition}[\cite{taylor}]\label{definition 1}
{\rm 
A continuous function $\gamma: S^n\to \mathbb{R}_+$ satisfying 
$\Gamma_\gamma=\mbox{inv}(\mbox{graph}(\gamma))$ is called a \textit{convex integrand}.   
}
\end{definition}

\begin{figure}[hbtp]
\begin{center}
\includegraphics[width=6cm]{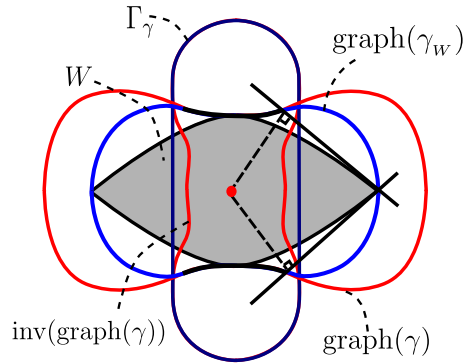}
\caption{$\gamma_{{}_W}(\theta)\le \gamma(\theta)$ for any $\theta\in S^n$.}
\label{figure 3}
\end{center}
\end{figure}

The following has been known.   
\begin{proposition}[\cite{taylor, nishimurasakemi2}]\label{proposition 1}
For any $\gamma_1, \gamma_2: S^n\to \mathbb{R}_+$, the following holds:  
\[
\Gamma_{\gamma_1}=\Gamma_{\gamma_2} 
%
\Leftrightarrow 
\mathcal{W}_{\gamma_1}= \mathcal{W}_{\gamma_2}.   
\]
\end{proposition}  
\noindent 
Proposition \ref{proposition 1} implies the following:   
\begin{proposition}\label{proposition 2}
Let $W$ be a convex body in $\mathbb{R}^{n+1}$ such that the origin of $\mathbb{R}^{n+1}$ is contained in {\rm int}$(W)$.    
Then, the following holds for any two $\gamma_1, \gamma_2\in C^0_W(S^n, \mathbb{R}_+)$: 
\[
\Gamma_{\gamma_1}=\Gamma_{\gamma_2}.   
\] 
\end{proposition}  
\noindent 
By Proposition \ref{proposition 2}, the following definition is well-defined:  
\begin{definition}
[\cite{taylor}]
\label{definition 2}
{\rm 
Let $W$ be a convex body in $\mathbb{R}^{n+1}$ such that the origin of $\mathbb{R}^{n+1}$ is contained in int$(W)$.    
Define the unique function $\gamma_{{}_W}: S^n\to \mathbb{R}_+$ as follows:   
\[
\mbox{graph}(\gamma_{{}_W})=\mbox{inv}(\Gamma_\gamma), 
\]
where $\gamma$ is a function of $C^0_W(S^n, \mathbb{R}_+)$.   
The function 
$\gamma_{{}_W}$ is called the \textit{convex integrand of} $W$.   
}
\end{definition} 
By the construction of convex integrand of $W$ and Proposition \ref{proposition 2}, 
the following holds:   
\begin{proposition}[\cite{taylor}]\label{proposition 3}
Let $W$ be a convex body containing the origin as an interior point.   
Then, the following holds for any $\theta\in S^n$ and any 
$\gamma\in C^0_W(S^n, \mathbb{R}_+)$.  
\[
\gamma_{{}_W}(\theta)\le \gamma(\theta).
\]  
\end{proposition}
\noindent 
Figure \ref{figure 3} illustrates Proposition \ref{proposition 3}. 
\subsection{Convex geometry in $S^{n+1}$}\label{convex in sphere}
\par 
For any point $P\in S^{n+1}$, let $H(P)$ be the closed hemisphere centered at $P$;  
namely, $H(P)$ is the set consisting of $Q\in S^{n+1}$ satisfying $P\cdot Q\ge 0$, 
where the dot in the center stands for the scalar product of two vectors 
$P, Q\in \mathbb{R}^{n+2}$.     
\begin{definition}[\cite{nishimurasakemi2}]\label{definition 3}
\rm{Let $\widetilde{W}$ be a subset of $S^{n+1}$.     
Suppose that there exists a point $P\in S^{n+1}$ such that 
$\widetilde{W}\cap H(P)=\emptyset$.      Then, $\widetilde{W}$ is said to be 
{\it hemispherical}.      
}
\end{definition}
For any non-empty subset $\widetilde{W}\subset S^{n+1}$, 
the {\it spherical polar set of $\widetilde{W}$}, denoted by 
$\widetilde{W}^\circ$, is defined as follows: 
\[
\widetilde{W}^\circ = \bigcap_{P\in \widetilde{W}}H(P).
\]   
\par 
\begin{lemma}[\cite{nishimurasakemi2}] 
For any hemispherical finite subset $\widetilde{X}=\{P_1, \ldots, P_k\}\subset S^{n+1}$, 
the following holds:
$$
\left\{\left.
\frac{\sum_{i=1}^k t_iP_i}{||\sum_{i=1}^kt_iP_i||}\;\right|\; 
P_i\in \widetilde{X},\; \sum_{i=1}^kt_i=1,\; t_i\ge 0
\right\}^\circ 
= 
H(P_1)\cap \cdots \cap H(P_k).
$$
\label{lemma 2.1}
\end{lemma}
Lemma \ref{lemma 2.1} is called {\it Maehara's lemma}.  
\par 
Let $P, Q$ be two points of $S^{n+1}$ such that $(1-t)P+tQ$ is not the zero vector 
for any $t\in [0,1]$.      Then, 
the following arc is denoted by $PQ$:  
\[
PQ=\left\{\left.\frac{(1-t)P+tQ}{||(1-t)P+tQ||}\in S^{n+1}\; \right|\; 0\le t\le 1\right\}.  
\]
\begin{definition}[\cite{hannishimura}] 
\label{definition 2.2}
{\rm 
Let $\widetilde{W}\subset S^{n+1}$ be a hemispherical subset.   
\begin{enumerate}
\item Suppose that 
$PQ\subset \widetilde{W}$ for any $P, Q\in \widetilde{W}$.     Then,  
$\widetilde{W}$ is said to be {\it spherical convex}.   
\item Suppose that $\widetilde{W}$ is closed, spherical convex and has an interior point.   
Then,  
$\widetilde{W}$ is said to be a {\it spherical convex body}.   
\end{enumerate}
}
\end{definition}  
%
\begin{definition}[\cite{hannishimura}]
{\rm 
Let $P$ be a point of $S^{n+1}$.   
\begin{enumerate}
\item A spherical convex body 
$\widetilde{W}$ such that $\widetilde{W}\cap H(-P)=\emptyset$ 
and $P\in \mbox{\rm int}(\widetilde{W})$ are satisfied is called 
a {\it spherical Wulff shape} relative to $P$.     
\item Let $\widetilde{W}$ be a spherical Wulff shape relative to $P$.    
Then, the set $\widetilde{W}^\circ$ 
is called the {\it spherical dual Wulff shape} of $\widetilde{W}$.   
\end{enumerate}
}
\end{definition}
\begin{definition}[\cite{nishimurasakemi2}]\label{definition 8}
{\rm Let $\widetilde{W}$ be a hemispherical subset of $S^{n+1}$.     
Then, the following set, denoted by 
$\mbox{\rm s-conv}(\widetilde{W})$, is called the {\it spherical convex hull of} 
$\widetilde{W}$.   
\[
\mbox{\rm s-conv}(\widetilde{W})= 
\left\{\left.
\frac{\sum_{i=1}^k t_iP_i}{||\sum_{i=1}^kt_iP_i||}\;\right|\; 
P_i\in \widetilde{W},\; \sum_{i=1}^kt_i=1,\; t_i\ge 0, k\in \mathbb{N}
\right\}.
\] 
}
\end{definition}
\begin{lemma}[\cite{nishimurasakemi2}]\label{lemma 2.2}
Let $\widetilde{W}_{1}, \widetilde{W}_{2}$ be non-empty subsets of $S^{n+1}$.    
Suppose that the inclusion $\widetilde{W}_{1}\subset \widetilde{W}_{2}$ holds.    
Then, the inclusion $\widetilde{W}_{2}^{\circ}\subset \widetilde{W}_{1}^{\circ}$ holds.
\end{lemma}
\begin{proposition}[\cite{nishimurasakemi2}]\label{proposition 7}
For any non-empty closed hemispherical subset 
$\widetilde{W} \subset S^{n+1}$, 
the equality $\mbox{ \rm s-conv}(\widetilde{W})= \left(
\mbox{ \rm s-conv}\left( \widetilde{W} \right) \right)^{\circ\circ }$ holds.
\end{proposition}
\subsection{Construction of Wulff shapes by using spherical polar sets}\label{spherical construction}
\par
\begin{figure}[hbtp]
\begin{center}
\includegraphics[width=6cm]{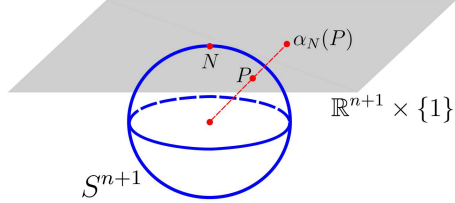}
\caption{The central projection $\alpha_N$.}
\label{figure 4}
\end{center}
\end{figure}

 
Let $Id: \mathbb{R}^{n+1}\to \mathbb{R}^{n+1}\times \{1\}\subset \mathbb{R}^{n+2}$ be the map defined by 
$Id(x)=(x,1)$.      
Denote the north-pole $(0, \ldots, 0, 1)\in \mathbb{R}^{n+2}$ by $N$.  
The set $S^{n+1}-H(-N)$ is denoted by $S_{N,+}^{n+1}$.    
Let $\alpha_N: S_{N,+}^{n+1}\to \mathbb{R}^{n+1}\times \{1\}$ be the central projection relative to 
$N$; 
namely, $\alpha_N$ is defined as follows where 
$P=(P_1, \ldots, P_{n+1}, P_{n+2})\in S_{N, +}^{n+1}$ (see Figure \ref{figure 4}):   
\[
\alpha_N\left(P_1, \ldots, P_{n+1}, P_{n+2}\right)
=
\left(\frac{P_1}{P_{n+2}}, \ldots, \frac{P_{n+1}}{P_{n+2}}\right).    
\]
\par 
Next, we consider the mapping $\Psi_N:S^{n+1}-\{\pm N\}\to {\color{black}S_{N, +}^{n+1}}$ 
defined by 
\[
\Psi_N(P)=\frac{1}{\sqrt{1-(N\cdot P)^2}}(N-(N\cdot P)P).    
\]
The mapping $\Psi_N$ was firstly introduced in \cite{nishimura}.    
It has been used for many purposes, 
for instance for the study of singularities of spherical pedal curves (\cite{nishimura, nishimura2}), 
for the study of pedal unfoldings of pedal curves (\cite{nishimura3}),  
for the study of hedgehogs (\cite{nishimurasakemi}) and 
for the study of a geometric model of crystal growth in the plane (\cite{kagatsumenishimura}).  
The hyperbolic version of $\Psi_N$ is also useful (see \cite{izumiyatari}).  
The mapping  $\Psi_N$  has the following 
intriguing properties:   
\begin{enumerate}
\item For any $P\in S^{n+1}-\{\pm N\}$, the equality $P\cdot \Psi_N(P)=0$ holds.    
\item For any $P\in S^{n+1}-\{\pm N\}$, 
the property $\Psi_N(P)\in \mathbb{R}N+\mathbb{R}P$ holds.      
\item For any $P\in S^{n+1}-\{\pm N\}$, the property $N\cdot \Psi_N(P)>0$ holds.       
\item The restriction $\Psi_N|_{S^{n+1}_{N,+}-\{N\}}: S^{n+1}_{N,+}-\{N\}\to S^{n+1}_{N,+}-\{N\}$ 
is a $C^\infty$ diffeomorphism.   
\end{enumerate}
Moreover, it is easily seen that by using $\Psi_N$, the inversion 
$\mbox{inv}: \mathbb{R}^{n+1}-\{0\}\to \mathbb{R}^{n+1}-\{0\}$ can be characterized as follows:  
\begin{proposition}\label{proposition 4}
\[
\mbox{\rm inv}=Id^{-1}\circ \alpha_N\circ \Psi_N\circ \alpha_N^{-1}\circ Id.   
\] 
\end{proposition}
\begin{proposition}[\cite{nishimurasakemi2}]\label{proposition 5}
Let $\gamma: S^{n}\to \mathbb{R}_+$ be a continuous function.    
Then, $\mathcal{W}_\gamma$ is characterized as follows:   
\[
\mathcal{W}_\gamma = 
Id^{-1}\circ \alpha_N\left(\left(\Psi_N\circ \alpha_N^{-1}\circ 
Id\left(\mbox{\rm graph}(\gamma)\right)\right)^\circ\right).   
\]
\end{proposition}
\begin{proposition}[\cite{nishimurasakemi2}]\label{proposition 6}
For any Wulff shape $\mathcal{W}_\gamma$, 
the following hold:   
\begin{enumerate}
\item The following set, too, is a Wulff shape.    
\[
Id^{-1}\circ \alpha_N\left(\left(\alpha_N^{-1}\circ 
Id\left(\mathcal{W}_\gamma\right)\right)^\circ\right).      
\] 
\item The graph of the convex integrand of $\mathcal{W}_\gamma$ is as follows.  
\[
\mbox{\rm inv}\left(
\partial \left(Id^{-1}\circ \alpha_N
\left(
\left(
\alpha_N^{-1}\circ Id(\mathcal{W}_\gamma) 
\right)^\circ 
\right)
\right)
\right).   
\]  
\item The graph of the convex integrand of 
$Id^{-1}\circ \alpha_N\left(\left(\alpha_N^{-1}\circ 
Id\left(\mathcal{W}_\gamma\right)\right)^\circ\right)$ is as follows, 
where $\partial\mathcal{W}_\gamma$ stands for the boundary of $\mathcal{W}_\gamma$.  
\[
\mbox{\rm inv}\left(
\partial 
\mathcal{W}_\gamma 
\right).   
\]  
\end{enumerate}
\end{proposition}
The assertions (2), (3) of Proposition \ref{proposition 6} has been implicitely 
proved in \cite{nishimurasakemi2}.   
\begin{definition}[\cite{nishimurasakemi2}]
{\rm 
For any Wulff shape $\mathcal{W}_\gamma$, the Wulff shape 
\[
Id^{-1}\circ \alpha_N\left(\left(\alpha_N^{-1}\circ 
Id\left(\mathcal{W}_\gamma\right)\right)^\circ\right)
\]  
given in Proposition \ref{proposition 6} 
is called the {\it dual Wulff shape} of 
$\mathcal{W}_\gamma$.   
}
\end{definition}
\section{Proof of Theorem \ref{theorem 1}}\label{section 3}
\subsection{Proof of the \lq\lq if\rq\rq part}
\indent    
In this subsection, we show that $W$ is strictly convex under the assumption that 
its convex integrand $\gamma_{{}_W}$ is 
of class $C^1$.   
Recall that $N$ (resp., $\alpha_N$) is the north-pole $(0, \ldots, 0,1)$ of $S^{n+1}$ (resp., the central projection relative to $N$).  
Set $\widetilde{W}=\alpha_N^{-1}\circ Id(W)$.   
By Proposition \ref{proposition 5}, we have the following:   
\[
W=Id^{-1}\circ \alpha_N\left(\left(
\Psi_N\circ \alpha_N^{-1}\circ Id
\left(\mbox{\rm graph}(\gamma_{{}_W})\right)\right)^\circ \right).
\]
Thus, we have the following:   
\[
\widetilde{W}=
\left(
\Psi_N\circ \alpha_N^{-1}\circ Id
\left(\mbox{\rm graph}(\gamma_{{}_W})\right)\right)^\circ.  
\]
By Proposition \ref{proposition 4}, we have the following:   
\[
\widetilde{W}=
\left(\alpha_N^{-1}\circ Id\circ \mbox{\rm inv}
\left(\mbox{\rm graph}(\gamma_{{}_W})\right)\right)^\circ.   
\]
By the definition of convex integrand and Proposition \ref{proposition 7}, 
the following holds:   
\[
\partial \widetilde{W}^\circ = 
\alpha_N^{-1}\circ Id\circ \mbox{\rm inv}
\left(\mbox{\rm graph}(\gamma_{{}_W})\right),   
\]
where $\partial \widetilde{W^\circ}$ stands for the boundary of $\widetilde{W^\circ}$.   
Thus, by the assumption of $\gamma_{{}_W}$, the following holds:   
\begin{lemma}\label{lemma 3.1}
The boundary of $\widetilde{W}^\circ$ is $C^1$ diffeomorphic to $S^n$. 
\end{lemma}
\par 
\smallskip 
Suppose that $W$ is not strictly convex.   
Then, since $W$ is convex, there exist two distinct points 
$x_0, x_1\in \partial W$ such that 
\[
\left\{
(1-t)x_0+tx_1\; |\; 0\le t\le 1
\right\}
\]
is contained in $\partial W$.    
For any $i\in \{0,1\}$, set $P_i=\alpha_N^{-1}(x_i)$.     
Then, it follows that 
\[
P_0P_1\subset \partial \widetilde{W}. 
\]
For any $\widetilde{P}\in S^{n+1}$, let 
$\partial H(\widetilde{P})$ be the boundary of $H(\widetilde{P})$.   
\begin{lemma}\label{lemma 3.2}
Under the above situation, the following holds:  
\[
\partial H(P_0)\cap \partial H(P_1)\cap \mbox{\rm int}(H(N))\subset \widetilde{W}^\circ.
\] 
\end{lemma}
\proof 
Let $P$ be a point of 
$\partial H(P_0)\cap \partial H(P_1)\cap \mbox{\rm int}(H(N))$.    
Suppose that $P\not\in \widetilde{W}^\circ$.    
Then, there exists a point $Q_0\in \widetilde{W}$ such that 
$P\cdot Q_0<0$.    
\par 
On the other hand, we have that 
$P\cdot P_i=0$ for any $i\in \{0,1\}$.    
Thus, it follows that 
\[
P\cdot 
\frac{(1-t)P_0+tP_1}{|| (1-t)P_0+tP_1||}
=0
\] 
for any $t\in [0,1]$.    
Since $P_0P_1\subset \partial \widetilde{W}$, 
$P\in \mbox{\rm int}\left(H(N)\right)$ and 
$\widetilde{W}$ is a spherical Wulff shape relative to $N$, it follows 
that $P\cdot Q\ge 0$ for any $Q\in \widetilde{W}$.   
\par 
Thus, we have a contradiction.   
\hfill {$\Box$}\\
\par 
\medskip 
By Lemma \ref{lemma 3.2}, 
for any $P\in \partial H(P_0)\cap \partial H(P_1)\cap \mbox{\rm int}(H(N))$, 
we have at least two great hyperspheres 
$\partial H(P_0), \partial H(P_1)$ 
which may be candidates for tangent great hyperspheres to 
$\widetilde{W}^\circ$ at $P$.   
This contradicts Lemma \ref{lemma 3.1}
\hfill {$\Box$}\\
\subsection{Proof of the \lq\lq only if\rq\rq part}
In this subsection, we show that $\gamma_{{}_W}$ is 
of class $C^1$ under the assumption that $W$ is strictly convex. 
We use the same notations given in Subsection 3.1.    
\begin{lemma}\label{lemma 3.3}
For any point $Q\in \partial \widetilde{W}^{\circ}$, 
there exists the unique point $P_Q\in \partial \widetilde{W}$
such that $Q$ is a point of $\partial H(P_Q)$.
\end{lemma}
\proof  
Suppose that 
there exists a point $Q\in \partial \widetilde{W}^{\circ}$ such that 
$Q\in \partial H(P_0)\cap \partial H(P_1)$, 
where $P_0, P_1$ are some distinct points of $\partial \widetilde{W}$.
Then, it follows that $P_i\cdot Q=0$ for any $i\in \{0,1\}$.    
This implies that $P\cdot Q= 0$ 
for any point $P \in P_1P_2$. 
Then, for any $\varepsilon >0$ 
and any point $P\in P_1P_2$, there exists a point $\widetilde{P}$
such that two inequalities $||P-\widetilde{P}||<\varepsilon$ and 
$\widetilde{P}\cdot Q<0$ are satisfied.  
Since $Q\in \widetilde{W}^\circ$, it follows that 
$\widetilde{P}\not \in \widetilde{W}$ although 
$P$ belongs to $\widetilde{W}$.   
Hence, we have that the arc $P_1P_2$ is contained in $\partial \widetilde{W}$.  
This contradicts the assumption that $W$ is strictly convex.   
\hfill {$\Box$}\\
\par 
\smallskip 
Secondly, we show the following lemma:  
\begin{lemma}\label{lemma 3.4}
The convex integrand $\gamma_{{}_W}: S^n\to \mathbb{R}_+$ 
is differentiable at any $\theta\in S^n$.    
\end{lemma}
\proof 
Let $Q$ be the point of $\partial \widetilde{W}^\circ\subset S^{n+1}$ such that the following 
is satisfied.   
\[
\mbox{\rm inv}\circ Id^{-1}\circ \alpha_N(Q)
= 
\left(\theta, \gamma_{{}_W}(\theta)\right).   
\]
Let $U$ be a sufficiently small neighbourhood 
of $Q$ in $S^{n+1}$.    
Since $\partial H(P_Q)$ is a great hypersphere and 
$\alpha_N$ is the central projection relative to $N$, 
we may assume that there exists an affine transformation 
$H: \mathbb{R}^{n+1}\to \mathbb{R}^{n+1}$ such that   
$H\circ Id^{-1}\circ \alpha_N(Q)$ is the origin of $\mathbb{R}^{n+1}$ and  
$H\circ Id^{-1}\circ \alpha_N(\partial H(P_Q)\cap U)$ (denoted by $V$)  
is contained in $\mathbb{R}^n\times \{0\}$, 
where $P_Q$ is the unique point of $\partial \widetilde{W}$ obtained in 
Lemma \ref{lemma 3.3}.
Let $V_1$ be a sufficiently small neighborhood of the origin in $V$.  
Then, by replacing $H=(h_1, \ldots, h_n, h_{n+1})$ with 
$(h_1, \ldots, h_n, -h_{n+1})$ if necessary, 
we may assume that 
there exists a continuous function $f_n: V_1\to [0, \infty)$ such that 
$f_n(0, \ldots, 0)=0$ and the graph of $f_n$
is an open subset of 
$H\circ Id^{-1}\circ \alpha_N(\partial \widetilde{W}^\circ\cap U)$.   
\par 
\smallskip 
We first suppose that $n=1$.   
Let $a$ be a positive real number such that 
$\{x\; |\; |x|< a\}\subset V_1 $.   
For the $a$, define $I(a)$ as follows: 
\[
I(a)=\{\lambda\in \mathbb{R}\; |\; 
\exists x\in (-a, 0)\cup (0, a) \mbox{ such that } 
f_1(x)=\lambda x\}.
\] 
The set $I(a)$ has the following properties:   
\begin{claim}\label{claim 3.1}
\begin{enumerate}
\item $0<a_1<a_2  \Rightarrow  0\le \sup I(a_1)\le \sup I(a_2)$.   
\item $a_2<a_1<0 \Rightarrow \inf I(a_2)\le \inf I(a_1)\le 0$. 
\item $\lim_{a\to 0}\sup I(a)=0$.   
\item $\lim_{a\to 0}\inf I(a)=0$.
\end{enumerate}
\end{claim}
\proof     
By definition, (1) and (2) are clear.   
\par 
We show (3).   
By (1), the following holds:   
\[
\lim_{a\to 0}\sup I(a)\ge 0.   
\]
Suppose that there exists a positive real number $\lambda_1$ such that 
$\lim_{a\to 0}\sup I(a)=\lambda_1$.   
Then, since $\widetilde{W}^\circ$ is spherical convex, 
the following holds for any $x$ such that $0<x<a$.   
\[
f_1(x)\ge \lambda_1x.   
\]
Since $\lambda_1>0$, the above inequality implies 
that there exists a point $P\in \partial \widetilde{W}$ $(P\ne P_Q)$ 
such that $Q\in \partial H(P)$.   
This contradicts Lemma \ref{lemma 3.3}.    
Therefore, we have $\lim_{a\to 0}\sup I(a)=0$.    
\par 
(4) may be proved similarly as (3).
\hfill {$\Box$}\\
\par 
\smallskip 
By Claim \ref{claim 3.1}, we have the following 
\begin{eqnarray*}
0 & \le & \lim_{x\to +0}\frac{f_1(x)}{x}\le \lim_{a\to 0}\sup I(a)=0, \\ 
0 & \ge & \lim_{x\to -0}\frac{f_1(x)}{x}\ge \lim_{a\to 0}\inf I(a)=0.  
\end{eqnarray*}
Therefore, $\gamma_{{}_W}$ must be differentiable at $\theta$.   
\par
\begin{figure}[hbtp]
\begin{center}
\includegraphics[width=6cm]{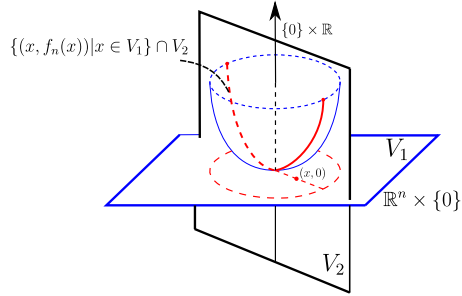}
\caption{$\{(x, f_n(x))\; |\; x\in V_1\}\cap V_2$ 
may be regarded as the graph of $f_1$.}
\label{figure 5}
\end{center}
\end{figure}

 
\medskip 
Next, we give a proof for general $n$.   
Let $x$ be a point of $V_1-\{0\}$.    
Set $V_2=\{0\}\times \mathbb{R}+\mathbb{R}(x,0)\subset \mathbb{R}^n\times \mathbb{R}$.  
Then, $V_2$ is a $2$-dimensional real vector space and 
the intersection 
\[
\{(x, f_n(x))\; |\; x\in V_1\}\cap V_2
\]
may be regarded as the graph of $f_1$ in the case $n=1$ (see Figure \ref{figure 5}).  
Let $\{r_i\}_{i=1, 2, \ldots}\subset \mathbb{R}_+$ be a sequence such that 
$\lim_{i\to \infty}r_i=0$.      
Set $x_i=r_ix$.     Then, by the proof in the case $n=1$, 
we have the following:  
\[
\lim_{i\to \infty}\frac{f_n(x_i)}{||x_i||}=0.   
\]    
Therefore, even for general $n$, $\gamma_{{}_W}$ must be differentiable at $\theta$. 
\hfill {$\Box$}\\
\par 
\medskip 
Let $Q$ be a point of 
$\partial \widetilde{W}^\circ$.     
Let $\{Q_i\}_{i=1, 2, \ldots}\subset \partial \widetilde{W}^\circ$ 
be a sequence such that $\lim_{i\to \infty}Q_i=Q$.   
By Lemma \ref{lemma 3.3}, for any $Q_i$ 
there exists the unique point $P_{Q_i}\in \partial \widetilde{W}$ 
such that 
$Q_i\in H(P_{Q_i})$.  
Set $P=P_Q$ and $P_i=P_{Q_i}$.          
By Lemma \ref{lemma 3.4}, in order to show that 
$\gamma_{{}_W}: S^n\to \mathbb{R}_+$ is of class $C^1$, 
it is sufficient to show the following:  
\[
\lim_{i\to \infty}h(H(P), H(P_i))=0, 
\] 
where $h: \mathcal{H}(S^{n+1})\times \mathcal{H}(S^{n+1})\to \mathbb{R}$ 
is the Pompeiu-Hausdorff metric.   
Suppose that there exists a positive real number 
$\varepsilon>0$ such that for any $m\in \mathbb{N}$ there exists 
an integer $i>m$ such that 
$h(H(P), H(P_i))>\varepsilon$.   
Since 
$Q\in H(P)$ and $Q_i\in H(P_i)$, by the definition of 
Pompeiu-Hausdorff metric (for the definition of Pompeiu-Hausdorff metric, 
see for instance \cite{barnsley, falconer}), 
it follows that 
 there exists a positive real number 
$\varepsilon>0$ such that for any $m\in \mathbb{N}$ there exists 
an integer $i>m$ such that 
$d(Q, Q_i)>\varepsilon$, 
where $d: \mathbb{R}^{n+2}\times \mathbb{R}^{n+2}\to \mathbb{R}$ is the Euclidean 
metric.    
This contradicts $\lim_{i\to \infty}Q_i=Q$.        
Therefore, we have $\lim_{i\to \infty}h(H(P), H(P_i))=0$. 
\hfill {$\Box$}
\section{Applications of Theorem \ref{theorem 1}}\label{section 4}
\begin{figure}[hbtp]
\begin{center}
\includegraphics[width=6cm]{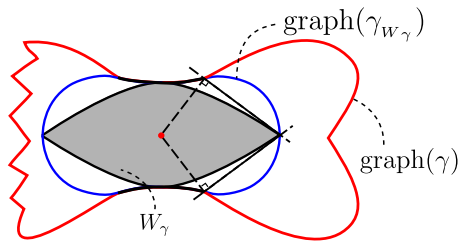}
\caption{A strictly convex $\mathcal{W}_{\gamma}$ having non smooth $\gamma$.}
\label{figure 6}
\end{center}
\end{figure}

\indent 
Since the boundary of the convex hull of a $C^1$ closed submanifold 
is a $C^1$ closed submanifold 
(for instance, see \cite{robertsonromerofuster, zakalyukin}), 
as a corollary of Theorem \ref{theorem 1}, we have the following:  
\begin{corollary}\label{corollary 1}
Let $\gamma: S^n\to \mathbb{R}_+$ be a function of class $C^1$.    
Then, $\mathcal{W}_\gamma$ is strictly convex.    
\end{corollary} 
\noindent 
In particular, we have the following:   
\begin{corollary}[\cite{nishimurasakemi2}, Theorem 1.3]\label{corollary 2}
Let $\gamma: S^n\to \mathbb{R}_+$ be a function of class $C^1$.    
Then, $\mathcal{W}_\gamma$ is never a polytope.    
\end{corollary}
\noindent 
On the other hand, Figure \ref{figure 6} shows 
that the converse of Corollary \ref{corollary 1} does not hold in general.   
\par 
Combining Theorem \ref{theorem 1} and Proposition \ref{proposition 6} 
yields the following:   
\begin{corollary}\label{corollary 3}
A Wulff shape in $\mathbb{R}^{n+1}$ is strictly convex 
if and only if the boundary of its dual Wulff shape 
is $C^1$ diffeomorphic to $S^n$.    
\end{corollary}  
\noindent 
In particular, we have the following:   
\begin{corollary}\label{corollary 4}
A Wulff shape in $\mathbb{R}^{n+1}$ is strictly convex and 
its boundary is $C^1$ diffeomorphic to $S^n$ 
if and only if its dual Wulff shape is strictly convex and the boundary of it 
is $C^1$ diffeomorphic to $S^n$.    
\end{corollary}  
\noindent 
It is interesting to compare Corollary \ref{corollary 4} and the following:   
\begin{proposition}[\cite{nishimurasakemi2}]
A Wulff shape in $\mathbb{R}^{n+1}$ is a polytope 
if and only if its dual Wulff shape is a polytope.  
\end{proposition}
Finally, we give an application of Theorem \ref{theorem 1} from the view point of 
pedal.   
\begin{definition}\label{definition 9}
{\rm Let $p$ (resp., $F: S^n\to \mathbb{R}^{n+1}$) be a point of 
$\mathbb{R}^{n+1}$ (resp., a $C^1$ embedding).  
Then, the {\it pedal of $F(S^n)$ relative to $p$} is the mapping 
$G: S^n\to \mathbb{R}^{n+1}$ which maps $\theta\in S^n$ to 
the nearest point in the tangent hyperplane to $F(S^n)$ at 
$F(\theta)$ from the given point $p$.  
}
\end{definition}
Let $W$ be a Wulff shape in $\mathbb{R}^{n+1}$.    
Suppose that $\partial W$ 
is $C^1$ diffeomorphic to 
$S^n$.   Then, $\partial W$
may be regarded as the graph of 
a certain $C^1$ embedding $F: S^n\to \mathbb{R}^{n+1}$, and    
$\gamma_{{}_W}$ is exactly the pedal of $\partial W$ relative to the origin.  
Theorem \ref{theorem 1} gives a sufficient condition for 
the pedal of $\partial W$ relative to the origin to be smooth:   
\begin{corollary}
Suppose that a Wulff shape $W$ in $\mathbb{R}^{n+1}$ is strictly convex and 
its boundary is $C^1$ diffeomorphic to $S^n$.    
Then, the pedal of $\partial W$ relative to the origin is of class $C^1$.    
\end{corollary} 
\section*{Acknowledgements}
The authors are grateful to Frank Morgan 
for helpful comments.   
This work is partially supported by JSPS KAKENHI 
Grant Number 26610035.


\begin{thebibliography}{99}          
{\color{black}
\bibitem{barnsley}M.~Barnsley, 
\textit{Fractals Everywhere Second edition}, Morgan Kaufmann Pub., San Fransisco, 1993.  
\bibitem{falconer}K.~Falconer, 
\textit{Fractal Geometry--Mathematical Foundations and applications  Second edition}, 
John Wiley \& Sons Ltd., Chichester, West Sussex, 2003.  
\bibitem{giga}Y.~Giga, 
\textit{Surface Evolution Equations}, Monographs of Mathematics, {\bf 99}, Springer, 2006.   
\bibitem{hannishimura}H.~Han and T.~Nishimura, 
\textit{The spherical dual transform is an isometry for spherical Wulff shapes}, 
preprint (available from arXiv:1504.02845 [math.MG]).  
\bibitem{izumiyatari}S.~Izumiya and F.~Tari, 
\textit{Projections of hypersurfaces in the hyperbolic space 
to hyperhorospheres and hyperplanes}, {Rev.~Mat.~Iberoam.} {\bf 24}(2008), 895--920.      
}
\bibitem{kagatsumenishimura} D.~Kagatsume and T.~Nishimura, \textit{Aperture of plane curves}, {J.~Singul.}, {\bf 12} (2015),80--91.  
{\color{black}
{\color{black}
}
\bibitem{morgan} F.~Morgan, 
\emph{The cone over the Clifford torus in $\mathbb{R}^4$ is F-minimizing}, 
Math. Ann., {\bf 289} (1991), no. 2, 341--354. 
\bibitem{nishimura}T.~Nishimura, 
\textit{Normal forms for singularities of pedal curves produced by non-singular 
dual curve germs in $S^n$},  
{Geom Dedicata} {\bf 133}(2008), 59--66.     
\bibitem{nishimura2}T.~Nishimura, 
\textit{Singularities of pedal curves produced by singular dual curve
germs in $S^n$}, Demonstratio Math., {\bf 43}(2010), 447--459.   
\bibitem{nishimura3}T.~Nishimura, 
\textit{Singularities of one-parameter pedal unfoldings of spherical pedal curves}, 
J.~Singul., {\bf 2}(2010), 160--169.   
\bibitem{nishimurasakemi}T.~Nishimura and Y.~Sakemi, 
{\it View from inside}, Hokkaido~Math.~J., {\bf 40}(2011), 361--373.        
}
\bibitem{nishimurasakemi2} T.\ Nishimura and Y.\ Sakemi, 
\textit{Topological aspect of Wulff shapes}, 
J. Math. Soc. Japan, {\bf 66} (2014), 89--109.
\bibitem{crystalbook}A.~Pimpinelli and J.~Villain, 
{\it Physics of Crystal Growth}, Monographs and Texts in Statistical Physics, 
Cambridge University Press, Cambridge New York, 1998. 
\bibitem{robertsonromerofuster}S.~A.~Robertson and M.~C.~Romero-Fuster, 
{\it The convex hull of a hypersurface}, 
Proc. London Math. Soc., {\bf 50}(1985), 370--384.   
\bibitem{schneider}R.~Schneider, 
\textit{Convex Bodies: The Brunn-Minkowski Theory 2nd edition}, Encyclopedia of Mathematics and its Applications  
{\bf 44}, Cambridge University Press, Cambridge, 2013.   
\bibitem{taylor}J.~E.~Taylor, \textit{Crystalline variational problems}, Bull. Amer. Math. Soc., 
{\bf 84}(1978), 568--588.     
\bibitem{taylor2}J.~E.~Taylor, J.~W.~Cahn and C.~A.~Handwerker, 
\textit{Geometric models of crystal growth}, 
Acta Metallurgica et Materialia, {\bf 40}(1992), 1443--1474.  
\bibitem{wulff}G.~Wulff, 
\textit{Zur frage der geschwindindigkeit 
des wachstrums und der aufl\"osung der krystallflachen}, 
Z. Kristallographine und Mineralogie, {\bf 34}(1901), 449--530.
\bibitem{zakalyukin}V.~M.~Zakalyukin, 
{\it Singularities of convex hulls of smooth maniifolds}, 
Functional Anal. Appl., {\bf 11}(1977), 225--227(1978).   
\end{thebibliography}
\end{document}